\documentclass[11pt,leqno]{article}
\usepackage{amsmath,amssymb,latexsym}
\usepackage{amscd}
\usepackage{mathrsfs}
\usepackage{amsthm}
\usepackage{amsfonts}

\newcommand{\R}{\mathbb R}
\newcommand{\N}{\mathbb N}

\newcommand{\C}{\mathbb C}

\newcommand{\To}{\longrightarrow}

\newcommand{\id}{\operatorname{id}}

\newcommand{\cP}{\mathcal P}
\newcommand{\cM}{\mathcal M}
\newcommand{\cN}{\mathcal N}
\newcommand{\im}{\operatorname{im}}
\newcommand{\alg}{\operatorname{alg}}

\let\epsilon\varepsilon
\let\phi\varphi

\theoremstyle{plain}
\newtheorem{Th}{Theorem}[section]
\newtheorem{Cor}[Th]{Corollary}
\newtheorem{Lem}[Th]{Lemma}
\newtheorem{Prop}[Th]{Proposition}
\newtheorem{OPr}[Th]{Open Problem}

\theoremstyle{definition}
\newtheorem{Def}[Th]{Definition}

\theoremstyle{remark}

\newtheorem{Rem}[Th]{Remark}

\title{{\sc On the algebra of smooth operators}}
\author{{\sc Tomasz Cia{\'s}}}
\date{}

\begin{document}

\maketitle

\begin{abstract}
Let $s$ be the space of rapidly decreasing sequences. 
We give the spectral representation of normal elements in the Fr\'echet algebra $L(s',s)$ of the so-called smooth operators. 
We also characterize closed commutative ${}^*$-subalgebras of $L(s',s)$ and establish a H\"older continuous functional calculus in this algebra. 
The key tool is the property $(DN)$ of $s$.
\end{abstract}


\section{Introduction}
The space $s$ of rapidly decreasing sequences plays significant role in the structure theory of nuclear Fr\'echet spaces. 
One of the most explicit example of this is provided by the K\=omura-K\=omura theorem which implies that a Fr\'echet space is nuclear if and only 
if it is isomorphic to some closed subspace of $s^\N$ (see e.g. \cite[Cor. 29.9]{MeV}). 
The space $s$ has also many interesting representations. For instance, it is isomorphic as a Fr\'echet space to 
the Schwartz space $\mathcal S(\R^n)$ of rapidly decreasing smooth functions, the space $\mathcal D(K)$ of test functions with support in 
a compact set $K\subset\R^n$ such that 
$\operatorname{int}(K)\neq\emptyset$, the space $C^\infty(M)$ of smooth functions on a compact smooth manifold $M$, the space $C^\infty[0,1]$ of smooth
functions on the interval $[0,1]$. Finally, the space $s$ and all of the spaces above are Fr\'echet commutative algebras with 
the pointwise multiplication. However, these algebras do not have to be isomorphic as algebras 
(for instance, $s$ and $C^\infty[0,1]$ with the pointwise multiplication are not isomorphic as algebras).

A natural candidate for ''the noncommutative $s$'' is the algebra $L(s',s)$ of the so-called smooth operators, where multiplication is just
the composition of operators (let us note that $s\subseteq s'$ continuously). It appears in $K$-theory for Fr\'echet algebras
(\cite[Def. 2.1]{Phill}, \cite[Ex. 2.12]{BhIO}, \cite{GlockLkamp}) and in $C^*$-dynamical systems (\cite[Ex. 2.6]{ElNatNest}). 
The algebra $L(s',s)$ is also an example of a dense smooth subalgebra of a $C^*$-algebra 
(precisely, $L(s',s)$ is a dense subalgebra of the $C^*$-algebra $K(l_2)$ of compact operators on $l_2$) 
which is especially important in the noncommutative geometry 
(see \cite{BhIO}, \cite{BlCun}, \cite[p. 23, 183-184]{Con}). From the philosophical point of view $C^*$-algebras just corespond to analogues of topological 
spaces whereas some of their dense smooth subalgebras play the role of smooth structures.
    
Representations of $s$ may lead up to representations of the algebra $L(s',s)$. Many of them are collected in \cite[Th. 2.1]{Dom}. For example,
$L(s',s)$ is isomorphic as a Fr\'echet ${}^*$-algebra to the following ${}^*$-algebras of continuous linear operators with the appropriate 
multiplication and involution: $L(\mathcal S'(\R^n), \mathcal S(\R^n))$, $L(\mathcal E'(M),C^\infty(M))$, 
$L(\mathcal E'[0,1],C^\infty[0,1])$, where $\mathcal S'(\R^n)$ is the space of tempered distributions, $M\subset \R^n$ is a compact smooth manifold,
$\mathcal E'(M)$ is the space of distributions on $M$ and $\mathcal E'[0,1]$ is the space of distributions with the support in $[0,1]$.
It is also worth to mention two extra representations of $L(s',s)$: the algebra of rapidly decreasing matrices 
\[\mathcal K:=\{(\xi_{j,k})_{j,k\in\N}: \sup_{j,k\in\N}|\xi_{j,k}|j^qk^q<\infty \text{ for all }q\in\N_0\}\]
with matrix multiplication and conjugation of the transpose as an involution (see e.g. \cite[p. 238]{Con}, \cite[Def. 2.1]{Phill}), 
and also the algebra $\mathcal S(\R^2)$ equipped with the Volterra convolution $(f\cdot g)(x,y):=\int_\R f(x,z)g(z,y)\mathrm{d}z$ and the involution 
$f^*(x,y):=\overline{f(y,x)}$ (see e.g. \cite[Ex. 2.12]{BhIO}).

The purpose of this paper is to present some spectral, algebra and functional calculus properties of the algebra of smooth operators. The results
are derived from the basic theory of nuclear Fr\'echet spaces and the theory of bounded operators on a separable Hilbert space.
The heart of the paper is the theorem on the spectral representation of normal elements in $L(s',s)$ (Theorem \ref{th_spectral}). In the proof we use 
the fact that the operator norm $||\cdot||_{l_2\to l_2}$ is a dominating norm on $L(s',s)$ (Proposition \ref{prop_l_2_DN}). 
As a by-product we obtain a kind of spectral description of normal elements of $L(s',s)$ among those of $K(l_2)$ (Corollary \ref{cor_characterization_of_L(s',s)}). 
Next, we characterize closed commutative ${}^*$-subalgebras of $L(s',s)$. We prove that every such a subalgebra is generated by a single operator 
and also by its spectral projections (Theorem \ref{th_commutative_subalg}), and, moreover, that it is a K\"othe algebra with pointwise multiplication. 
To do this, we show that every closed commutative ${}^*$-subalgebra of $L(s',s)$ 
has some canonical Schauder basis (Lemma \ref{lem_minimal_projections}).  
Finally, we establish a H\"older-continuous functional calculus in  $L(s',s)$ (Corollary \ref{cor_holder})
and we prove the functional calculus theorem for normal elements in this algebra  
(Theorem \ref{th_functional_calculus}). 

By a Fr\'echet space we mean a complete metrizable locally convex space. A Fr\'echet algebra is a Fr\'echet space which is an algebra with a continuous 
multiplication. A Fr\'echet ${}^*$-algebra is a Fr\'echet algebra with an involution. 

We use the standard notation and terminology. All the notions from Functional Analysis are explained in \cite{MeV} and those from topological algebras
in  \cite{Fra} or \cite{Zel}.

\section{Preliminaries}\label{preliminaries}
Throughout the paper, $\N$ will denote the set of natural numbers $\{1,2,\ldots\}$ and $\N_0:=\N\cup\{0\}$.

By projection on $l_2$ we always mean a continuous orthogonal (self-adjoint) projection. 

We define \emph{the space of rapidly decreasing sequences} as a Fr\'echet space
\[s:=\bigg\{\xi=(\xi_j)_{j\in\N}\in\C^\N:|\xi|_q:=\bigg(\sum_{j=1}^\infty|\xi_j|^2j^{2q}\bigg)^{\frac{1}{2}}<\infty\text{ for all }q\in\N_0\bigg\}\]
with the topology corresponding to the system $(|\cdot|_q)_{q\in\N_0}$ of norms. 
Its strong dual is isomorphic to \emph{the space of slowly increasing sequences} 
\[s':=\bigg\{\xi=(\xi_j)_{j\in\N}\in\C^\N:|\xi|_q^{'}:=\bigg(\sum_{j=1}^\infty|\xi_j|^2j^{-2q}\bigg)^{\frac{1}{2}}<\infty\text{ for some }q\in\N_0\bigg\}\]
equipped with the inductive limit topology given by the system $(|\cdot|_q^{'})_{q\in\N_0}$ of norms.

Every $\eta\in s'$ corresponds to the continuous functional $\xi\mapsto \langle \xi,\eta \rangle$ on $s$, where
\[\langle \xi,\eta\rangle:=\sum_{j=1}^\infty \xi_j\overline{\eta_j}.\]
Futhermore, by the Cauchy-Schwartz inequality we get
\[|\langle \xi,\eta\rangle|\leq |\xi|_q|\eta|'_q\] 
for all $q\in\N_0$, $\xi\in s$ and $\eta\in s'$ with $|\eta|'_q<\infty$.

For $1\leq p<\infty$ and the K\"othe matrix $(a_{j,q})_{j\in\N,q\in\N_0}$ we define the K\"othe space 
\[\lambda^p(a_{j,q}):=\bigg\{\xi=(\xi_j)_{j\in\N}\in\C^\N:|\xi|_{p,q}:=\bigg(\sum_{j=1}^\infty|\xi_j a_{j,q}|^p\bigg)^{\frac{1}{p}}<\infty\text{ for all }q\in\N_0\bigg\}\]
and for $p=\infty$
\[\lambda^\infty(a_{j,q}):=\bigg\{\xi=(\xi_j)_{j\in\N}\in\C^\N:|\xi|_{\infty,q}:=\sup_{j\in\N}|\xi_j| a_{j,q}<\infty\text{ for all }q\in\N_0\bigg\}\]
with the topology generated by norms $(|\cdot|_{p,q})_{q\in\N_0}$ (see e.g. \cite[Def. on p. 326]{MeV}). Please note that sometimes these spaces are 
Fr\'echet ${}^*$-algebras with the pointwise multiplication. 

Now, the space $s$ is just 
the K\"othe space $\lambda^2(j^q)$. Moreover, since the space $s$ is a nuclear Fr\'echet space, it is isomorphic to any K\"othe space 
$\lambda^p(j^q)$ for $1\leq p\leq\infty$ (see e.g. \cite[Prop. 28.16, Ex. 29.4 (1)]{MeV}).  
We use $l_2$-norms to simplify futher computations, for example we have $|\cdot|_0=|\cdot|_0^{'}=||\cdot||_{l_2}$. 
  
It is well known that the space $L(s',s)$ of continuous linear operators from $s'$ to $s$ with the fundamental system of norms 
$(||\cdot||_q)_{q\in\N_0}$,
\[||x||_q:=\sup_{|\xi|_q^{'}\leq 1}{|x\xi|_q}\]
is isomorphic to $s$ as Fr\'echet space. 
Moreover, $L(s',s)$ is isomorphic to $s\widehat\otimes s$, the completed tensor product of $s$ (see \cite[\S41.7 (5)]{Kot}).

Since the canonical inclusion $j:s\hookrightarrow s'$ is continuous, thus for $x,y\in L(s',s)$
\[x\cdot y:=x\circ j\circ y\]
is in $L(s',s)$ as well and with this operation $L(s',s)$ is a Fr\'echet algebra. 

The diagram
\[l_2\hookrightarrow s'\to s\hookrightarrow l_2\]
defines the canonical (continuous) embeding of the algebra $L(s',s)$ into the algebra $L(l_2)$ of continuous linear operators on the Hilbert space $l_2$. 
In fact, this inclusion acts into the space $K(l_2)$ of compact operators on $l_2$ and the sequence of singular numbers of elements in $L(s',s)$
belongs to $s$ (see \cite[Prop. 3.1, Cor. 3.2]{Dom}). 
Therefore, $L(s',s)$ can be regarded as some class of compact operators on $l_2$. 
Clearly, the multiplication in $L(s',s)$ coincides with the composition in $L(l_2)$, and further, $L(s',s)$ is invariant under the hilbertian involution 
$x\mapsto x^*$. 

To see this, let us consider the Fr\'echet ${}^*$-algebra of rapidly decreasing matrices 
\[\mathcal K:=\{\Xi=(\xi_{j,k})_{j,k\in\N}:|||\Xi|||_q:=\sup_{j,k\in\N}|\xi_{j,k}|j^qk^q<\infty \text{ for all }q\in\N_0\},\]
with the matrix multiplication, the involution defined by $((\xi_{j,k})_{j,k\in\N})^*:=(\overline{\xi_{k,j}})_{j,k\in\N}$ and with 
$(|||\cdot|||_q)_{q\in\N_0}$ as its fundamental sequence of norms. 
By \cite[Th. 2.1]{Dom}, $\Phi\colon L(s',s)\to \mathcal{K}$, $\Phi(x):=(\langle xe_k,e_j\rangle)_{j,k\in\N}$ is the algebra isomorphism
and we have 
\[\Phi(x)^*=\big(\overline{\langle xe_j,e_k\rangle}\big)_{j,k\in\N}=(\langle x^*e_k,e_j\rangle)_{j,k\in\N}.\]
Hence, $x^*=\Phi^{-1}(\Phi(x)^*)\in L(s',s)$ and $\Phi$ is even a ${}^*$-isomorphism. Clearly, for every matrix $\Xi\in\mathcal{K}$ and $q\in\N_0$, $|||\Xi^*|||_q=|||\Xi|||_q$, thus the hilbertian involution
is continuous on $L(s',s)$.

The Fr\'echet algebra $L(s',s)$ with the involution ${}^*$ is called \emph{the algebra of smooth operators}. We will also consider 
the algebra with unit
\[\widetilde{L(s',s)}:=\{x+\lambda\textbf{1}:x\in L(s',s),\lambda\in\C\},\]
where \textbf{1} is the identity operator on $l_2$. We endow the algebra $\widetilde{L(s',s)}$ with the product topology.



Now, we shall recall some basic spectral properties of the algebra $L(s',s)$. 
For the sake of convenience, we state the following definition.
\begin{Def}\label{def_eigenvalues}
We say that a sequence $(\lambda_n)_{n\in\N}\subset\C$ is a \emph{sequence of eigenvalues} of an infinite dimensional compact operator $x$ on $l_2$ if 
it satisfies the following conditions:
\begin{enumerate}
 \item[(i)] $\{\lambda_n\}_{n\in\N}$ is the set of eigenvalues of $x$ without zero;
 \item[(ii)] $|\lambda_1|\geq|\lambda_2|\geq\ldots> 0$ and if two eigenvalues has the same absolute value then we can ordered them in an arbitrary way;  
 \item[(iii)] the number of occurrences of the eigenvalue $\lambda_n$ is equal to its geometric multiplicity
(i.e., the dimension of the space $\ker(\lambda_n\textbf{1}-x)$).
\end{enumerate}
\end{Def}

Proposition \ref{prop_4} is well known 
(see e.g. \cite{GlockLkamp}) and it is a simple consequence of Proposition \ref{prop_3}. 
Proofs of Propositions \ref{prop_3}, \ref{prop_5} one can find in \cite[Th. 3.3, Cor. 3.5]{Dom}.
         

\begin{Prop}\label{prop_3}
An operator in $\widetilde{L(s',s)}$ is invertible if and only if it is invertible in $L(l_2)$.     
\end{Prop}

\begin{Prop}\label{prop_4}
The algebra $\widetilde{L(s',s)}$ is a $Q$-algebra, i.e., the set of invertible elements is open. 
\end{Prop}

\begin{Prop}\label{prop_5}
The spectrum of $x$ in $L(s',s)$ equals the spectrum of $x$ in $L(l_2)$ and it consists of zero and the set of 
eigenvalues. If moreover, $x$ is infinite dimensional, then the sequence of eigenvalues of $x$ (see Definition \ref{def_eigenvalues}) belongs to $s$.    
\end{Prop}

The first part of the following proposition is also known (see e.g. \cite[Lemma 2.2]{Phill}). We give a simple proof that norms 
$||\cdot||_q$ are submultiplicative. 
  
\begin{Prop}\label{prop_6}
The algebra $L(s',s)$ is $m$-locally convex, i.e., it has a fundamental system of submultiplicative norms. 
Moreover, $||xy||_q\leq||x||_q||y||_q$ for every $q\in\N_0$.
\end{Prop}
\proof
Let $x,y\in L(s',s)$ and let $B_q$, $B_q'$ denote the closed unit ball for the norms $|\cdot|_q$, $|\cdot|_q'$, respectively. Clearly, 
$y(B_q')\subseteq||y||_q B_q$ and $B_q\subseteq B_q'$. Hence
\begin{align*}
||xy||_q&=\sup_{|\xi|_q'\leq1}|x(y(\xi))|_q=\sup_{\eta\in y(B_q')}|x(\eta)|_q\leq\sup_{\eta\in ||y||_q B_q}|x(\eta)|_q
=||y||_q\sup_{\eta\in B_q}|x(\eta)|_q\\
&\leq ||y||_q\sup_{\eta\in B'_q}|x(\eta)|_q=||x||_q||y||_q. 
\end{align*}
\qed

\section{Spectral representation of normal elements}

In this section we prove the following theorem on the spectral representation of normal elements in $L(s',s)$.
\begin{Th}\label{th_spectral}
Every infinite dimensional normal operator $x$ in $L(s',s)$ has the unique spectral representation $x=\sum_{n=1}^\infty\lambda_nP_n$, where 
$(\lambda_n)_{n\in\N}$ is a decreasing (according to the modulus) sequence in $s$ of nonzero pairwise different elements, 
$(P_n)_{n\in\N}$ is a sequence of nonzero pairwise orthogonal finite dimensional projections belonging to $L(s',s)$ 
(i.e., the canonical inclusion of $P_n$ into $L(l_2)$ is a projection on $l_2$) 
and the series converges absolutely in $L(s',s)$. 
Moreover, $(|\lambda_n|^\theta||P_n||_q)_{n\in\N}\in s$ for all $q\in\N_0$ and all $\theta\in(0,1]$. 
\end{Th}

Before we prove this result, we first need to do some preparation. 
Let us recall (see \cite[Def. on p. 359 and Lemma 29.10]{MeV}) that a Fr\'echet space $(X,(||\cdot||_q)_{q\in\N_0})$ has the 
\emph{property $(DN)$} if there is a continuous norm $||\cdot||$ on $X$ such that 
for any $q\in\N_0$ and $\theta\in(0,1)$ there is $r\in\N_0$ and $C>0$ such that for all $x\in X$
\[||x||_q\leq C||x||^{1-\theta}||x||_r^\theta.\]
The norm $||\cdot||$ is called a \emph{dominating norm}.

The following result is essentially due to K. Piszczek (see \cite[Th. 4]{Piszczek}).
\begin{Prop}\label{prop_l_2_DN}
The norm $||\cdot||_{l_2\to l_2}$ is a dominating norm on $L(s',s)$.
\end{Prop}
\proof
Clearly, $||\cdot||_{l_2\to l_2}=||\cdot||_0$. By \cite[Th. 4.3]{V1} (see the proof), the conclusion is equivalent to the following condition 
\[\forall q\in\N_0\phantom{i}\forall\theta>0\phantom{i}\exists r\in\N_0\phantom{i}\exists C>0\phantom{i}\forall h>0\phantom{i}\forall x\in L(s',s)\quad
||x||_q\leq C\bigg(h^\theta||x||_r+\frac{1}{h}||x||_0\bigg).\]
From H\"older's inequality, the norm $|\cdot|_0$ is a dominating norm on $s$, hence again by \cite[Th. 4.3]{V1} as above we get
\[\forall q\in\N_0\phantom{i}\forall\eta>0\phantom{i}\exists r\in\N_0\phantom{i}\exists D_0>0\phantom{i}\forall k>0\phantom{i}\forall \xi\in s\quad
|\xi|_q\leq D_0\bigg(k^\eta|\xi|_r+\frac{1}{k}|\xi|_0\bigg).\]
Now, by the bipolar theorem (see e.g. \cite[Th. 22.13]{MeV}) we obtain (follow the proof of \cite[Lemma 29.13]{MeV}) an equivalent condition
\begin{equation}\label{eq_DNsubset}
\forall q\in\N_0\phantom{i}\forall\eta>0\phantom{i}\exists r\in\N_0\phantom{i}\exists D>0\phantom{i}\forall k>0\quad
U_q^\circ\subset D\bigg(k^\eta U_r^\circ+\frac{1}{k}U_0^\circ\bigg),
\end{equation}
where $U_q:=\{\xi\in s:|\xi|_q\leq1\}$ and $U_q^\circ$ is its polar.
If $\theta>0$ and $h\in(0,1]$ are given, we define $\eta:=2\theta+1$, $k:=\sqrt{h}$. Since $k^{2\eta}\leq k^{\eta-1}$, we obtain
\begin{align*}
U_q^\circ\otimes U_q^\circ&:=\{x\otimes y:x,y\in U_q^\circ\}\subset D\bigg(k^\eta U_r^\circ+\frac{1}{k}U_0^\circ\bigg)\otimes D\bigg(k^\eta U_r^\circ+\frac{1}{k}U_0^\circ\bigg)\\
&\subset D^2\bigg(k^{2\eta}U_r^\circ\otimes U_r^\circ+2k^{\eta-1}U_r^\circ\otimes U_r^\circ+\frac{1}{k^2}U_0^\circ\otimes U_0^\circ\bigg)\\
&\subset 3D^2\bigg(k^{\eta-1}U_r^\circ\otimes U_r^\circ+\frac{1}{k^2}U_0^\circ\otimes U_0^\circ\bigg)
=3D^2\bigg(h^{\theta}U_r^\circ\otimes U_r^\circ+\frac{1}{h}U_0^\circ\otimes U_0^\circ\bigg).
\end{align*}
Since $r$ and $D$ in condition (\ref{eq_DNsubset}) can be choosen so that $q\leq r$ and $D\geq1$, thus for $h>1$ it holds
\[U_q^\circ\otimes U_q^\circ\subset U_r^\circ\otimes U_r^\circ
\subset 3D^2\bigg(h^\theta U_r^\circ\otimes U_r^\circ+\frac{1}{h}U_0^\circ\otimes U_0^\circ\bigg),\]
hence we have shown
\[\forall q\in\N_0\phantom{i}\forall\theta>0\phantom{i}\exists r\in\N_0\phantom{i}\exists C>0\phantom{i}\forall h>0\quad
U_q^\circ\otimes U_q^\circ\subset C\bigg(h^\theta U_r^\circ\otimes U_r^\circ+\frac{1}{h}U_0^\circ\otimes U_0^\circ\bigg).\]
Therefore, for all $x:=\sum_{j=1}^nx_j\otimes y_j\in s\otimes s$ we get
\begin{align*}
\sup\{|z(x)|:z\in U_q^\circ\otimes U_q^\circ\}
&\leq C\sup\bigg\{|z(x)|:z\in h^\theta U_r^\circ\otimes U_r^\circ+\frac{1}{h}U_0^\circ\otimes U_0^\circ\bigg\}\\
&=C\sup\bigg\{|(z'+z'')(x)|:z'\in h^\theta U_r^\circ\otimes U_r^\circ,z''\in\frac{1}{h}U_0^\circ\otimes U_0^\circ\bigg\}\\
&\leq C\sup\bigg\{|z'(x)|+|z''(x)|:z'\in h^\theta U_r^\circ\otimes U_r^\circ,z''\in\frac{1}{h}U_0^\circ\otimes U_0^\circ\bigg\}\\
&=C\bigg(h^\theta\sup\{|z(x)|:z\in U_r^\circ\otimes U_r^\circ\}
+\frac{1}{h}\sup\{|z(x)|:z\in U_0^\circ\otimes U_0^\circ\}\bigg).
\end{align*}
Let $\chi\colon s\otimes s\to L(s',s)$, $\chi\big(\sum_{j=1}^nx_j\otimes y_j\big)(z):=\sum_{j=1}^n z(y_j)x_j.$
We have for all $p\in\N_0$
\begin{align*}
\sup\bigg\{\bigg|z\bigg(\sum_{j=1}^nx_j\otimes y_j\bigg)\bigg|:z\in U_p^\circ\otimes U_p^\circ\bigg\}
&=\sup\bigg\{\bigg|\sum_{j=1}^nz_1(x_j) z_2(y_j)\bigg|:z_1,z_2\in U_p^\circ\bigg\}\\
&=\sup\bigg\{\bigg|z_1\bigg(\sum_{j=1}^n z_2(y_j)x_j\bigg)\bigg|:z_1,z_2\in U_p^\circ\bigg\}\\
&=\sup\bigg\{\bigg|\sum_{j=1}^n z(y_j)x_j\bigg|_p:z\in U_p^\circ\bigg\}\\
&=\sup\bigg\{\bigg|\chi\bigg(\sum_{j=1}^nx_j\otimes y_j\bigg)(z)\bigg|_p:z\in U_p^\circ\bigg\}\\
&=\bigg|\bigg|\chi\bigg(\sum_{j=1}^nx_j\otimes y_j\bigg)\bigg|\bigg|_p.
\end{align*}
Hence
\[\bigg|\bigg|\chi\bigg(\sum_{j=1}^n x_j\otimes y_j\bigg)\bigg|\bigg|_q\leq C\bigg(h^\theta\bigg|\bigg|\chi\bigg(\sum_{j=1}^nx_j\otimes y_j\bigg)\bigg|\bigg|_r
+\frac{1}{h}\bigg|\bigg|\chi\bigg(\sum_{j=1}^nx_j\otimes y_j\bigg)\bigg|\bigg|_0\bigg).\]

Finally, since the set $\{\chi\big(\sum_{j=1}^nx_j\otimes y_j\big):x_j,y_j\in s, n\in\N\}$ is dense in $L(s',s)$, thus
\[||x||_q\leq C\bigg(h^\theta||x||_r+\frac{1}{h}||x||_0\bigg)\]
for all $x\in L(s',s)$.
\qed

\begin{Lem}\label{lem_1}
Let $(E, (||\cdot||_q)_{q\in\N_0})$ be a Fr\'echet space with the property $(DN)$ and let $||\cdot||_p$ be a dominating norm. 
If $(x_n)_{n\in\N}\subset E$, $(\lambda_n)_{n\in\N}\subset\C$ satisfy conditions
\begin{enumerate}
 \item[(i)] $\sup_{n\in\N}{||x_n||_p}<\infty$,
 \item[(ii)] $\forall q\in\N_0\phantom{i}\sup_{n\in\N}{|\lambda_n|||x_n||_q}<\infty$, 
\end{enumerate}
then
\[\forall q\in\N_0\phantom{i}\forall\theta\in(0,1]\quad\sup_{n\in\N}{|\lambda_n|^\theta||x_n||_q}<\infty.\]
Moreover, for another sequence $(y_n)_{n\in\N}\subset E$ satisfying conditions (i) and (ii) we have
\[\forall q\in\N_0\phantom{i}\forall q'\in\N_0\phantom{i}\forall\theta\in(0,1]\quad\sup_{n\in\N}{|\lambda_n|^\theta||x_n||_q||y_n||_{q'}}<\infty.\]
\end{Lem}
\proof
Let us fix $q\in\N_0$ and $\theta\in(0,1)$. Since $||\cdot||_p$ is a dominating norm on $E$, thus for some $C>0$ and $r\in\N_0$
\begin{equation}\label{eq_DN_x_n}
||x_n||_q\leq C||x_n||_p^{1-\theta}||x_n||_r^\theta
\end{equation}
for all $n\in\N$. Let $C_1:=\sup_{n\in\N}||x_n||_p<\infty$, $C_2:=\sup_{n\in\N}|\lambda_n|||x_n||_q<\infty$. Then by (\ref{eq_DN_x_n}) 
\[|\lambda_n|^\theta||x_n||_q\leq C||x_n||_p^{1-\theta}(|\lambda_n|||x_n||_r)^\theta\leq CC_1^{1-\theta}C_2^\theta=:C_3,\]
where $C_3$ does not depend on $n$.

To prove the second assertion let us also fix $q'\in\N_0$ and let $(y_n)_{n\in\N}\subset E$ satisfy conditions (i) and (ii). We have 
\[|\lambda_n|^\theta||x_n||_q||y_n||_{q'}=(|\lambda_n|^{\frac{\theta}{2}}||x_n||_q)(|\lambda_n|^{\frac{\theta}{2}}||y_n||_{q'})\]
and from the first part of the proof, $\sup_{n\in\N}|\lambda_n|^{\frac{\theta}{2}}||x_n||_q<\infty$ and 
$\sup_{n\in\N}|\lambda_n|^{\frac{\theta}{2}}||y_n||_{q'}<\infty$, 
so we are done.  
\qed

\begin{Prop}\label{prop_P_n_in_L(s',s)}
Let $\cN$ be either a finite set or $\cN=\N$. 
If $(P_n)_{n\in\cN}$ is a sequence of pairwise orthogonal finite dimensional projections on $l_2$, $(\lambda_n)_{n\in\cN}\subset\C\setminus\{0\}$ and 
$x:=\sum_{n\in\cN}\lambda_nP_n\in L(s',s)$ (the series converges in norm $||\cdot||_{l_2\to l_2}$), then $(P_n)_{n\in\cN}\subset L(s',s)$. 
\end{Prop}
\proof
Since, $P_n=\frac{1}{\lambda_n}x\circ P_n$ thus $P_n\colon l_2\to s$, but $P_n=P_n\circ\frac{1}{\lambda_n}x$ so $P_n$ extends to $P_n\colon s'\to l_2$. 
Hence, $P_n=P_n\circ P_n\colon s'\to s$.
\qed

\begin{Lem}\label{lem_2}
Let $(\lambda_n)_{n\in\N}$ be a decreasing (according to the modulus) sequence of nonzero complex numbers and let $(P_n)_{n\in\N}$ be a sequence of 
nonzero pairwise orthogonal finite dimensional projections on $l_2$.
Let us also assume that the series $\sum_{n=1}^\infty\lambda_nP_n$ converges in the norm $||\cdot||_{l_2\to l_2}$ and its limit belongs to $L(s',s)$. Then 
$(\lambda_n)_{n\in\N}\in s$, $(P_n)_{n\in\N}\subset L(s',s)$ and the series converges absolutely in $L(s',s)$. 
Moreover, $(|\lambda_n|^\theta||P_n||_q)_{n\in\N}\in s$ for all $q\in\N_0$ and $\theta\in(0,1]$.
\end{Lem}
\proof
By Prop. \ref{prop_5}, the sequence of eigenvalues of the operator $x:=\sum_{n=1}^\infty\lambda_nP_n$ belongs to $s$. Clearly, $\lambda_n$ is an eigenvalue
of $\sum_{n=1}^\infty\lambda_nP_n$ and the number of its occurrences is less or equal to the geometric multiplicity so $(\lambda_n)_{n\in\N}$ is, 
likewise, in $s$.   

By Prop. \ref{prop_P_n_in_L(s',s)}, $P_n\in L(s',s)$. 
We will show that 
$(|\lambda_n|^\theta||P_n||_q)_{n\in\N}\in s$ for all $q\in\N_0$ and $\theta\in(0,1]$, which implies that the series $\sum_{n=1}^\infty\lambda_nP_n$ converges
absolutely in $L(s',s)$. Let us consider operator $T_x\colon L(l_2)\to L(s',s)$ which sends $z\in L(l_2)$ to the following composition (in $L(s',s)$):
\[s'\stackrel{x}{\to}s\hookrightarrow l_2\stackrel{z}{\to}l_2\hookrightarrow s'\stackrel{x}{\to}s.\]
By the closed graph theorem for Fr\'echet spaces (see e.g. \cite[Th. 24.31]{MeV}), $T_x$ is continuous and since the sequence of operators 
$(P_n)_{n\in\N}$ is bounded in $L(l_2)$, thus 
$(\lambda_n^2P_n)_{n\in\N}=(T_xP_n)_{n\in\N}$ is bounded in $L(s',s)$, i.e.,
\begin{equation}\label{eq_||P_n||_q}
 \sup_{n\in\N}|\lambda_n|^2||P_n||_q<\infty
\end{equation}
for all $q\in\N_0$. 

Let $(e_n)_{n\in\N}$ be the canonical orthonormal basis in $l_2$ and let $E_n\colon s'\to s$,
\[E_n\xi:=\xi_ne_n,\]
for $\xi=(\xi_n)_{n\in\N}\in s'$ and $n\in\N$. Clearly, $E_n$ are projections in $L(s',s)$. Moreover,
\begin{equation}\label{eq_||E_n||_q=n^2q}
||E_n||_q=\sup_{|\xi|'_q\leq1}|E_n\xi|_q=\sup_{|\xi|'_q\leq1}|\xi_ne_n|_q=\sup_{|\xi|'_q\leq1}|\xi_n|\cdot|e_n|_q=n^q\cdot n^q=n^{2q}. 
\end{equation}
Since $(\lambda_n)_{n\in\N}\in s$, thus
\begin{equation}\label{eq_||E_n||_q}
 \sup_{n\in\N}|\lambda_n|^2||E_n||_q<\infty
\end{equation}
for $q\in\N_0$

By Proposition \ref{prop_l_2_DN}, $||\cdot||_{l_2\to l_2}$ is a dominating norm on $L(s',s)$, and, of course, $||P_n||_{l_2\to l_2}=||E_n||_{l_2\to l_2}=1$
for $n\in\N$. Thus, from inequalities (\ref{eq_||P_n||_q}), (\ref{eq_||E_n||_q}), equality (\ref{eq_||E_n||_q=n^2q}) and by Lemma \ref{lem_1} 
(applied to sequences $(\lambda_n^2)_{n\in\N}$, $(P_n)_{n\in\N}$ and $(E_n)_{n\in\N}$) we get
\[\sup_{n\in\N}|\lambda_n|^{2\theta}||P_n||_qn^{2q'}=\sup_{n\in\N}|\lambda_n|^{2\theta}||P_n||_q||E_n||_{q'}<\infty\]
for all $\theta\in(0,1]$ and $q,q'\in\N_0$. Hence, $(|\lambda_n|^\theta||P_n||_q)\in s$ for all $q\in\N_0$ and $\theta\in(0,1]$. 
\qed

Now, it is not hard to prove the main theorem of this section. 

\emph{Proof of Theorem \ref{th_spectral}}.\phantom{i}
Let $x$ be a normal infinite dimensional operator in $L(s',s)$. Operator $x$ (as an operator on $l_2$) is compact (see \cite[Prop. 3.1]{Dom}), 
thus by the spectral theorem for normal compact operators (see e.g. \cite[Th. 7.6]{Conw}), $x=\sum_{n=1}^\infty\lambda_nP_n$, where $(\lambda_n)_{n\in\N}$ is a decreasing 
zero sequence of nonzero pairwise different elements, $(P_n)_{n\in\N}$ is a sequence of nonzero pairwise orthogonal finite dimensional projections 
and the series converges in norm $||\cdot||_{l_2\to l_2}$. Now, the conclusion follows by Lemma \ref{lem_2}.   
\qed

As a corollary, we get a characterization of normal operators in $L(s',s)$ among compact operators on $l_2$.  

\begin{Cor}\label{cor_characterization_of_L(s',s)}
Let $x$ be a compact infinite dimensional normal operator on $l_2$ with the spectral representation $x=\sum_{n=1}^\infty\lambda_n P_n$ (the series
converges in norm $||\cdot||_{l_2\to l_2}$). 
Then the following assertions are equivalent:
\begin{enumerate}
 \item[(i)] $x\in L(s',s)$ (as an operator on $l_2$);
 \item[(ii)] $P_n\in L(s',s)$ for $n\in\N$ and $(|\lambda_n|^\theta||P_n||_q)_{n\in\N}\in s$ for all $q\in\N_0$ and every $\theta\in(0,1]$;
 \item[(iii)] $P_n\in L(s',s)$ for $n\in\N$, $(\lambda_n)_{n\in\N}\in s$ and $\sup_{n\in\N}|\lambda_n|||P_n||_q<\infty$ for all $q\in\N_0$;
 \item[(iv)] $P_n\in L(s',s)$ for $n\in\N$ and $\sum_{n=1}^\infty|\lambda_n|||P_n||_q<\infty$ for all $q\in\N_0$.
\end{enumerate} 
Moreover, if $x=\sum_{n=1}^N\lambda_nP_n$ is a finite dimensional operator on $l_2$, 
then $x\in L(s',s)$ if and only if $P_n\in L(s',s)$ for $n=1,\ldots,N$. 
\end{Cor}

\proof
The implication (i)$\Rightarrow$(ii) follows directly by Theorem \ref{th_spectral}. The implications (ii)$\Rightarrow$(iii), (iv)$\Rightarrow$(i) are obvious.

(iii)$\Rightarrow$(iv) By Lemma \ref{lem_1}, $\sup_{n\in\N}|\lambda_n|^{\frac{1}{2}}||P_n||_q<\infty$, and, moreover, $\sum_{n=1}^\infty|\lambda_n|^{\frac{1}{2}}<\infty$ ($s\subset \bigcap_{p>0}l_p$). Thus
\[\sum_{n=1}^\infty|\lambda_n|||P_n||_q\leq \sup_{n\in\N}|\lambda_n|^{\frac{1}{2}}||P_n||_q\cdot\sum_{n=1}^\infty|\lambda_n|^{\frac{1}{2}}<\infty,\]
and the remaining implication (iii)$\Rightarrow$(iv) holds. 

The finite case is an immediate consequence of Proposition \ref{prop_P_n_in_L(s',s)}.
\qed

\section{Closed commutative ${}^*$-subalgebras of $L(s',s)$}

The aim of this section is to describe all closed commutative ${}^*$-subalgebras of $L(s',s)$ (see Theorem \ref{th_commutative_subalg}) and 
identifying maximal among them (see Theorem \ref{th_max_com}).


We will need the following lemma.
\begin{Lem}\label{lem_a_1_ldots_a_n_in_A}
Let A be a subalgebra of the algebra $\widetilde A$ over $\C$. Let $N\in\N$, $a_1,\ldots,a_N\in\widetilde A$, 
$\lambda_1,\ldots,\lambda_N\in\C$, $a_j\neq0$, $a_j^2=a_j$, $a_ja_k=0$ for $j\neq k$, $\lambda_j\neq 0$ and $\lambda_j\neq\lambda_k$ for $j\neq k$. Then $a_1,\ldots,a_N\in A$ whenever $\lambda_1a_1+\ldots+\lambda_Na_N\in A$.
\end{Lem}
\proof
We use induction with respect to $N$. The case $N=1$ is trivial. 

Let us assume that the conlusion holds for all $M<N$. 
Let $a:=\lambda_1a_1+\ldots+\lambda_Na_N\in A$. We have
\[\lambda_1^2a_1+\ldots+\lambda_N^2a_N=a^2\in A,\]
and, on the other hand,
\[\lambda_N\lambda_1a_1+\ldots+\lambda_N^2a_N=\lambda_Na\in A\]
so
\[(\lambda_1^2-\lambda_N\lambda_1)a_1+\ldots+(\lambda_{N-1}^2-\lambda_N\lambda_{N-1})a_{N-1}=a^2-\lambda_Na\in A.\]
Since $\lambda_j\neq0$ and $\lambda_j\neq\lambda_N$ for $j\in\{1,\ldots,N-1\}$, thus $\lambda_j^2-\lambda_N\lambda_j=\lambda_j(\lambda_j-\lambda_N)\neq0$
for $j\in\{1,\ldots,N-1\}$. If $\lambda_j^2-\lambda_N\lambda_j$ are pairwise different then, from the inductive assumption, $a_1,\ldots,a_{N-1}\in A$ so $a_N\in A$ as well.

Let us assume that it is not the case.
We define the equivalence relation $\mathcal R$ on the set $\{1,\ldots,N-1\}$ in the following way
\[j\mathcal R k\Leftrightarrow \lambda_j(\lambda_j-\lambda_N)=\lambda_k(\lambda_k-\lambda_N).\]
Let $I_1,\ldots,I_{N_1}$ denote the equivalence classes which contain not less than two elements and let $I_0:=\{i_1,\ldots,i_{N_0}\}$ be the set of indices
from the remaining equivalence classes. From our assumption, $I_1\neq\emptyset$. For $j\in\{1,\ldots,N_1\}$, $k\in I_j$ let 
\[\lambda_j':=\lambda_k(\lambda_k-\lambda_N)\] 
and let 
\[a_j':=\sum_{n\in I_j}a_n.\] 
We also define 
\[\lambda_{N_1+1}':=\lambda_{i_1}(\lambda_{i_1}-\lambda_N),\lambda_{N_1+2}':=\lambda_{i_2}(\lambda_{i_2}-\lambda_N),\ldots,\lambda_{N_1+N_0}':=\lambda_{i_{N_0}}(\lambda_{i_{N_0}}-\lambda_N)\] 
and 
\[a_{N_1+1}':=a_{i_1},a_{N_1+2}':=a_{i_2},\ldots,a_{N_1+N_0}':=a_{i_{N_0}}.\] 
Clearly, $1\leq N':=N_1+N_0<N$, $a_j'\neq0$, $a_j'^2=a_j'$, $a_j'a_k'=0$, $\lambda_j'\neq0$, 
$\lambda_j'\neq\lambda_k'$ for $j,k\in\{1,\ldots,N'\}$, $j\neq k$ and 
\[\lambda_1'a_1'+\ldots+\lambda_{N'}'a_{N'}'=a^2-\lambda_Na\in A.\]
From the inductive assumption, $a_1'\in A$, hence 
\[\sum_{n\in I_1}\lambda_na_n=\sum_{n\in I_1}a_n\cdot\sum_{n=1}^N\lambda_na_n=a_1'a\in A.\] 
Again, from the inductive assumption, $a_n\in A$ for $n\in I_1$, and, therefore, 
$\sum_{n\in\{1,\ldots,N\}\setminus I_1}\lambda_na_n\in A$. Once again, from the inductive assumption, $a_n\in A$ for $n\in\{1,\ldots,N\}\setminus I_1$. Thus
$a_1,\ldots,a_{N}\in A$ which completes the proof.
\qed

\begin{Prop}\label{prop_x_in_A_iff_P_n_in_A}
Let $A$ be a closed ${}^*$-subalgebra of $L(s',s)$ (not necessary commutative) and let $x$ be an infinite dimensional normal operator in $L(s',s)$ with the spectral 
representation $x=\sum_{n=1}^\infty\lambda_nP_n$. Then $x\in A$ if and only if $P_n\in A$ for all $n\in\N$.
\end{Prop}
\proof
By Theorem \ref{th_spectral}, if $P_n\in A$ for all $n\in\N$ then $x\in A$. 
To prove the converse let us assume that $x\in A$. Then 
$x^*=\sum_{n=1}^\infty\overline{\lambda_n}P_n\in A$ so $xx^*=\sum_{n=1}^\infty|\lambda_n|^2P_n\in A$. 
Let $N_0:=0$, $N_1:=\sup\{n\in\N:|\lambda_n|=|\lambda_1|\}$ and for $j=2,3,\ldots$, let $N_j:=\sup\{n\in\N:|\lambda_n|=|\lambda_{N_{j-1}+1}|\}$. 
Since $(|\lambda_n|)_{n\in\N}$ is a zero sequence thus $N_j<\infty$. We have   
\[y_k:=\bigg(\frac{xx^*}{|\lambda_1|^2}\bigg)^k=\sum_{n=1}^\infty\bigg(\frac{|\lambda_n|}{|\lambda_1|}\bigg)^{2k}P_n\in A\]
for all $k\in\N$. For $q$ and $k$ arbitrary, we get
\begin{align*}
||y_k-(P_1+\ldots+P_{N_1})||_q&=\bigg|\bigg|\sum_{n=1}^\infty\bigg(\frac{|\lambda_n|}{|\lambda_1|}\bigg)^{2k}P_n-(P_1+\ldots+P_{N_1}) \bigg|\bigg|_q
=\bigg|\bigg|\sum_{n=N_1+1}^\infty\bigg(\frac{|\lambda_n|}{|\lambda_1|}\bigg)^{2k}P_n\bigg|\bigg|_q\\
&\leq\sum_{n=N_1+1}^\infty\bigg(\frac{|\lambda_n|}{|\lambda_1|}\bigg)^{2k}||P_n||_q
\leq\frac{1}{|\lambda_1|}\bigg(\frac{|\lambda_{N_1+1}|}{|\lambda_1|}\bigg)^{2k-1}\sum_{n=N_1+1}^\infty|\lambda_n|||P_n||_q.
\end{align*}
By Theorem \ref{th_spectral}, $\sum_{n=N_1+1}^\infty|\lambda_n|||P_n||_q<\infty$, and, moreover, $\frac{|\lambda_{N_1+1}|}{|\lambda_1|}<1$, thus 
\[||y_k-(P_1+\ldots+P_{N_1})||_q\to0\] 
as $k\to\infty$. Therefore, since $A$ is closed,
we conlude that $P_1+\ldots+P_{N_1}\in A$. 
Consequently, 
\[\sum_{n=N_1+1}^\infty|\lambda_n|^2P_n=xx^*-|\lambda_1|^2(P_1+\ldots+P_{N_1})\in A\]
hence, proceeding by induction, we get that $P_{N_j+1}+\ldots P_{N_{j+1}}\in A$ for $j\in\N_0$ so
\[\sum_{n=N_j+1}^{N_{j+1}}\lambda_nP_n=(P_{N_j+1}+\ldots P_{N_{j+1}})x\in A.\]
Finally, by Lemma \ref{lem_a_1_ldots_a_n_in_A}, $P_n\in A$ for $n\in\N$
\qed

\begin{Prop}\label{prop_series_converges_in_norm}
For every othonormal system $(e_n)_{n\in\N}$ in $l_2$ and sequence $(\lambda_n)_{n\in\N}\in c_0$, the series $\sum_{n=1}^\infty\lambda_n\langle\cdot,e_n\rangle e_n$ converges in norm $||\cdot||_{l_2\to l_2}$.
\end{Prop}
\proof
This is a simple consequence of the Pythagorean theorem and the Bessel's inequality.
\qed

\begin{Lem}\label{lem_minimal_projections}
Let $A$ be a commutative subalgebra of $L(s',s)$. Let $\cP$ denote the set of nonzero projections belonging to $A$ and 
let $\cM$ be the set of minimal elements in $\cP$ with respect to the following order relation
\[\forall P,Q\in\cP \quad P\preceq Q\Leftrightarrow PQ=QP=P.\]
Then\\
(i) $\cM$ is at most countable family of pairwise orthogonal projections belonging to $L(s',s)$ such that 
\[\forall P\in\cP\phantom{i}\exists P'_1,\ldots,P'_m\in\cM\quad P=P'_1+\cdots+P'_m.\]
(ii) If $A$ is also a closed ${}^*$-subalgebra of $L(s',s)$, then $\cM$ is a Schauder basis in $A$. 
\end{Lem}
\proof
(i) By the definition
\[\cM=\{P\in\cP:\forall Q\in\cP\quad(Q\preceq P\Rightarrow Q=P)\}.\]
Firstly, we will show that
\begin{equation}\label{eq_forallP_existsQ_j}
\forall P\in\cP\phantom{i}\exists P'_1,\ldots,P'_m\in\cM\quad P=P'_1+\cdots+P'_m.
\end{equation}

Let us take $P\in\cP$. If $P\in\cM$, then we are done. Otherwise, there is $Q\in\cP$ such that $Q\preceq P$, $Q\neq P$. Of course, $P-Q\in\cP$. 
If $Q,P-Q\in\cM$, then $P=Q+(P-Q)$ is desired decomposition. Otherwise, we decompose $Q$ or $P-Q$ into smaller projections as it was done above for 
projection $P$. Since $P$ is finite dimensional thus after finitely many steps we finish our procedure.

Next, we shall prove that projections in $\cM$ are pairwise orthogonal. Let $P,Q\in\cM$, $P\neq Q$ and let us suppose, to derive a contradiction, that $PQ\neq0$. 
Since $A$ is commutative, thus
\[(PQ)^2=P^2Q^2=PQ\]
so $PQ\in\cP$. Moreover,
\[P(PQ)=P^2Q=PQ\]
so $PQ\preceq P$. Now, $PQ\neq P$ implies that $P\notin\cM$ and if $PQ=P$ then $Q\notin\cM$, which is a contradiction.

Finally, since projections in $\cM$ are pairwise orthogonal (as projections on $l_2$) thus $\cM$ is at most countable.

(ii) Let $x\in A$. If $x$ is finite dimensional and $\sum_{n=1}^N\mu_nQ_n$ is its spectral decomposistion, then from (i) and by Lemma \ref{lem_a_1_ldots_a_n_in_A}, 
$x$ is a linear combination of projections in $\cM$.  

Let us assume that $x$ is infinite dimensional and let $x=\sum_{n=1}^\infty\mu_nQ_n$ (the spectral representation of $x$). 
Since $A$ is a closed commutative ${}^*$-subalgebra of $L(s',s)$ thus, by Proposition \ref{prop_x_in_A_iff_P_n_in_A}, $Q_n\in A$ for $n\in\N$. Next, from (i), we have that
\[\forall n\in\N\phantom{i}\exists Q^{(n)}_1,\ldots,Q^{(n)}_{l_n}\in\cM\quad Q_n=\sum_{j=1}^{l_n}Q^{(n)}_j.\]
Hence 
\[x=\sum_{n=1}^\infty\sum_{j=1}^{l_n}\mu_nQ^{(n)}_j.\]
For $l_0=0$, $j=l_0+l_1+\ldots+l_{n-1}+k$, $1\leq k\leq l_n$ let $P_j:=Q^{(n)}_k$ and let $\lambda_j:=\mu_{n}$. 
Let us consider the series $\sum_{n=1}^\infty\lambda_nP_n$. Clearly, if the series converges in $L(s',s)$ then its limit is $x$.  
To prove this we shall firstly show that the series $\sum_{n=1}^\infty\lambda_nP_n$ converges in the norm $||\cdot||_{l_2\to l_2}$. 

Since $P_n$ is a (orthogonal) projection of finite dimension $d_n$ thus $P_n=\sum_{j=1}^{d_n}\langle\cdot,e_j^{(n)}\rangle e_j^{(n)}$ for every 
orthonormal basis $\big(e_j^{(n)}\big)_{j=1}^{d_n}$ of the image of 
$P_n$. For $d_0=0$, $j=d_0+d_1+\ldots+d_{n-1}+k$, $1\leq k\leq d_{n}$ let $e_j:=e^{(n)}_k$ and let $\lambda_j':=\lambda_{n}$. 
By Proposition \ref{prop_series_converges_in_norm}, the series
$\sum_{j=1}^\infty\lambda_j'\langle\cdot,e_j\rangle e_j$ converges in the norm $||\cdot||_{l_2\to l_2}$. 
Hence $\sum_{n=1}^\infty\lambda_nP_n$ converges in the norm $||\cdot||_{l_2\to l_2}$ 
because $(\sum_{n=1}^N\lambda_nP_n)_{N\in\N}$ is a subsequence of the sequence of partial sums of the series 
$\sum_{j=1}^\infty\lambda_j'\langle\cdot,e_j\rangle e_j$.

Now, by Lemma \ref{lem_2}, $x=\sum_{n=1}^\infty\lambda_nP_n$ and the series converges absolutely in $L(s',s)$. This show that every operator in $A$
is represented by the absolutely convergent series $\sum_{n=1}^\infty\lambda_n''P_n''$ with $P_n''\in\cM$. To prove the uniquness of this representation
let us assume that $\sum_{n=1}^\infty\lambda_n''P_n''=0$. Then 
\[\lambda_m''P_m''=P_m''\sum_{n=1}^\infty\lambda_n''P_n''=0\]
so $\lambda_m''=0$ for $m\in\N$. This show that the sequence of coefficients is unique, hence $\cM$ is a Schauder basis in $A$.  
\qed

For a closed commutative ${}^*$-subalgebra $A$ of the algebra $L(s',s)$ the Schauder basis $\cM$ from the preceeding lemma will be called 
\emph{the canonical Schauder basis} (of $A$).

For a subset $Z$ of $L(s',s)$ we will denote by $\alg(Z)$ the closed ${}^*$-subalgebra of $L(s',s)$ generated by $Z$. 
If $A$ is a closed ${}^*$-subalgebra of $L(s',s)$, then $\widehat A$ denotes the set of nonzero ${}^*$-multiplicative functionals on $A$. 

\begin{Cor}\label{cor_coefficient_functionals}
The set $\widehat A$ of nonzero ${}^*$-multiplicative functionals on a closed commutative ${}^*$-subalgebra $A$ of $L(s',s)$ is exactly the set 
of coefficient functionals with respect to the canonical Schauder basis of $A$. 
\end{Cor}
\proof
Clearly, every coefficient functional is ${}^*$-multiplicative. Conversly, if $\phi$ is a nonzero ${}^*$-multiplicative functional on $A$ and 
$\{P_n\}_{n\in\N}$ is the canonical Schauder basis then 
$\phi(P_n)=\phi(P_n^2)=(\phi(P_n))^2$, 
thus $\phi(P_n)=0$ or $\phi(P_n)=1$. Let us suppose that $\phi(P_n)=\phi(P_m)=1$ for $n\neq m$. Then
\begin{align*}
2&=\phi(P_n)+\phi(P_m)=\phi(P_n+P_m)=\phi((P_n+P_m)^2)=(\phi(P_n+P_m))^2\\
&=(\phi(P_n)+\phi(P_m))^2=4,
\end{align*}
a contradiction. Hence, there is at most one $n\in\N$ such that $\phi(P_n)=1$. If $\phi(P_n)=0$ for all $n\in\N$ then, since $\{P_n\}_{n\in\N}$ is
a basis, $\phi=0$, a contradiction. Thus, there is exactly one $n\in\N$ such that $\phi(P_n)=1$, and 
$\phi(P_m)=0$ for $m\neq n$, i.e., $\phi$, is a coefficient functional.
\qed

\begin{Prop}\label{prop_alg=lin}
If $\{P_n\}_{n\in\cN}$ is a family of pairwise orthogonal projections belonging to $L(s',s)$, then
\[\alg(\{P_n\}_{n\in\cN})=\overline{\operatorname{lin}}(\{P_n\}_{n\in\cN})\]
and it is a commutative ${}^*$-algebra.
\end{Prop}
\proof
Clearly, $\overline{\operatorname{lin}}(\{P_n\}_{n\in\cN})\subseteq\alg(\{P_n\}_{n\in\cN})$ and $\operatorname{lin}(\{P_n\}_{n\in\cN})$ 
is a commutative ${}^*$-algebra. By the continuity of the algebra multiplication and the
hilbertian involution, $\overline{\operatorname{lin}}(\{P_n\}_{n\in\cN})$ is a commutative ${}^*$-algebra as well.
Hence, $\overline{\operatorname{lin}}(\{P_n\}_{n\in\cN})=\alg(\{P_n\}_{n\in\cN})$.
\qed

\begin{Prop}\label{prop_P_n_basic_sequnce}
Every sequence $\{P_n\}_{n\in\cN}\subset L(s',s)$ of nonzero pairwise orthogonal projections is a basic sequence in $L(s',s)$, i.e., 
it is a (canonical) Schauder basis of the Fr\'echet space (${}^*$-algebra) $\overline{\operatorname{lin}}(\{P_n\}_{n\in\cN})$.
\end{Prop}
\proof
Let $\cM$ be the canonical Schauder basis of $A:=\alg(\{P_n\}_{n\in\cN})$ which consists of all projections which are minimal with respect to the order 
relation described in Lemma \ref{lem_minimal_projections}.
If we show that $\{P_n\}_{n\in\cN}=\cM$, then, by Proposition \ref{prop_alg=lin}, we get the conclusion. 

Let $n\in\cN$ be fixed and let us assume that $Q\preceq P_n$ for some nonzero projection $Q\in A$, i.e., 
$QP_n=Q$. Since $A=\overline{\operatorname{lin}}(\{P_n\}_{n\in\cN})$ thus 
\[Q=\lim_{j\to\infty}\sum_{k=1}^{M_j}\lambda_k^{(j)}P_k\]
for some $M_j\in\N$ and $\lambda_k^{(j)}\in\C$. From the continuity of the algebra and scalar multiplication, we get
\begin{align*}
Q=QP_n&=\bigg(\lim_{j\to\infty}\sum_{k=1}^{M_j}\lambda_k^{(j)}P_k\bigg)P_n=\lim_{j\to\infty}\bigg(\sum_{k=1}^{M_j}\lambda_k^{(j)}P_kP_n\bigg)
=\lim_{j\to\infty}(\lambda_n^{(j)}P_n)\\
&=(\lim_{j\to\infty}\lambda_n^{(j)})P_n=\lambda_nP_n,
\end{align*}
where $\lambda_n:=\lim_{j\to\infty}\lambda_n^{(j)}\in\C$. Since $Q$ is a nonzero projection, thus 
$\lambda_n=1$ and $Q=P_n$. Hence $\{P_n\}_{n\in\cN}\subseteq\cM$.

Now, let us suppose that there is a projection $Q$ in $\cM\setminus\{P_n\}_{n\in\cN}$. We have already proved that $\{P_n\}_{n\in\cN}\subseteq\cM$, hence by Lemma \ref{lem_minimal_projections} (i), $Qx=0$ for all 
$x\in\operatorname{lin}(\{P_n\}_{n\in\cN})$. By continuity of the multiplication, 
$Qx=0$ for every $x\in\overline{\operatorname{lin}}(\{P_n\}_{n\in\cN})=A$. In particular, $Q=Q^2=0$, a contradiction. Hence, $\{P_n\}_{n\in\cN}=\cM$. 
\qed

Closed commutative ${}^*$-subalgebras of $L(s',s)$ are quite simple, all of them are generated by a single operator and also by its spectral projections. 
From nuclearity we get also the following sequence space representations.

\begin{Th}\label{th_commutative_subalg}
Let $A$ be a closed commutative infinite dimensional ${}^*$-subalgebra of $L(s',s)$. Then $A$ has a (canonical) Schauder basis $\{P_n\}_{n\in\N}$ 
consisting of pairwise othogonal finite dimensional minimal projections (see Lemma \ref{lem_minimal_projections}) such that
\[A=\alg(\{P_n\}_{n\in\N})\cong\lambda^1(||P_n||_q)=\lambda^\infty(||P_n||_q)\] 
as Fr\'echet ${}^*$-algebras. Moreover, there is an operator $x\in A$ with the spectral representation 
$x=\sum_{n=1}^\infty\lambda_nP_n$ such that $A=\alg(x)$.  
\end{Th}
\proof
By Lemma \ref{lem_minimal_projections}, $A$ has a Schauder basis with the desired properties. 
By Proposition \ref{prop_alg=lin}, $A=\overline{\operatorname{lin}}(\{P_n\}_{n\in\N})=\alg(\{P_n\}_{n\in\N})$ 
and since $A$ is a nuclear Fr\'echet space with the Schauder basis $\{P_n\}_{n\in\N}$, thus
(see e.g. \cite[Cor. 28.13, Prop. 28.16]{MeV}) 
\[A\cong\lambda^1(||P_n||_q)=\lambda^\infty(||P_n||_q)\]
as Fr\'echet spaces. Since on linear span of $\{P_n\}_{n\in\N}$ the multiplication (involution) corresponds to the pointwise multiplication 
(conjugation) in $\lambda^1(||P_n||_q)$ thus the isomorphism is also ${}^*$-algebra isomorphism where the K\"othe 
space is equipped with the pointwise multiplication.

Now, we shall show that there is a decreasing sequence $(\lambda_n)_{n\in\N}$ of positive numbers such that the series 
$\sum_{n=1}^\infty\lambda_nP_n$ is absolutely convergent in $L(s',s)$. 
To do so, let us choose a sequence $(C_q)_{q\in\N}$ such that $C_q\geq\max_{1\leq n\leq q}||P_n||_q$. 
Clearly, $\frac{C_q}{||P_n||_q}\geq1$ for $q\geq n$ so 
\[\inf_{q\in\N}\frac{C_q}{||P_n||_q}\geq\min\bigg\{\frac{C_1}{||P_n||_1},\ldots,\frac{C_{n-1}}{||P_n||_{n-1}},1\bigg\}>0\] 
for $n\in\N$. Let $\lambda_1:=1$ and let 
\[\lambda_n:=\min\bigg\{\frac{1}{n^2}\inf_{q\in\N}\frac{C_q}{||P_n||_q},\frac{\lambda_{n-1}}{2}\bigg\}.\]
Then $\lambda_n>0$, the sequence $(\lambda_n)_{n\in\N}$ is strictly decreasing and
\[\sum_{n=1}^\infty\lambda_n||P_n||_q\leq\sum_{n=1}^\infty\frac{1}{n^2}\inf_{r\in\N}\frac{C_r}{||P_n||_r}||P_n||_q
\leq C_q\sum_{n=1}^\infty\frac{1}{n^2}<\infty.\]

Consequently, $x:=\sum_{n=1}^\infty\lambda_nP_n\in L(s',s)$ and this series is the spectral representation of $x$. 
Moreover, since $P_n\in A$ for $n\in\N$ and $A$ is closed, thus $x\in A$.
Finally, the equality $\alg(x)=\alg(\{P_n\}_{n\in\N})$ is a consequence of Proposition \ref{prop_x_in_A_iff_P_n_in_A}. 
\qed 

A commutative closed ${}^*$-subalgebra $A$ of the algebra $L(s',s)$ is said to be \emph{maximal commutative} if it is not contain in any larger closed 
commutative ${}^*$-subalgebra of $L(s',s)$. We say that the sequence $\{P_n\}_{n\in\N}$ of nonzero pairwise orthogonal projections in $L(s',s)$ 
is \emph{complete} if there is no nonzero projection $P$ in $L(s',s)$ such that $P_nP=0$ for every $n\in\N$.
For a subset $Z$ of the algebra $L(s',s)$, the set $Z':=\{x\in L(s',s):xy=yx \text{ for all }y\in Z\}$ is called the \emph{commutant} 
of $Z$. 

\begin{Prop}\label{prop_Z'}
For every subset $Z$ of $L(s',s)$, the comutant $Z'$ is a closed ${}^*$-subalgebra of $L(s',s)$. 
\end{Prop}
\proof
Clearly, if $x,y$ commute with every $z\in Z$ then $\lambda x$, $x+y$, $xy$ and $x^*$ commute as well. Hence, from the continuity of the algebra operations 
and the involution, $Z'$ is a closed ${}^*$-subalgebra of $L(s',s)$. 
\qed

\begin{Th}\label{th_max_com}
For a closed commutative ${}^*$-subalgebra $A$ of $L(s',s)$ the following assertions are equivalent:
\begin{enumerate}
 \item[(i)] $A$ is maximal commutative;
 \item[(ii)] The canonical Schauder basis $\{P_n\}_{n\in\N}$ of $A$ is a complete sequence of pairwise orthogonal one-dimensional projections 
 belonging to $L(s',s)$;
 \item[(iii)] $A=A'$.
\end{enumerate} 
\end{Th}
\proof
(i)$\Rightarrow$(ii). 
Let us suppose that  for some $m\in\N$ the projection $P_m$ is not one dimensional. Then there are two nonzero pairwise orthogonal 
projections $Q_1,Q_2\in L(s',s)$ such that $P_m=Q_1+Q_2$. 
By Proposition \ref{prop_alg=lin}, $\overline{\operatorname{lin}}(\{P_n:n\neq m\}\cup\{Q_1,Q_2\})$ is a closed commutative ${}^*$-subalgebra of 
$L(s',s)$, and, clearly,
\[A=\overline{\operatorname{lin}}(\{P_n\}_{n\in\N})\subseteq \overline{\operatorname{lin}}(\{P_n:n\neq m\}\cup\{Q_1,Q_2\}).\] 
By Proposition \ref{prop_P_n_basic_sequnce}, $\{P_n\}_{n\in\N}$ is the canonical Schauder basis of $A$ and $\{P_n:n\neq m\}\cup\{Q_1,Q_2\}$ is the canonical Schauder basis of 
$\overline{\operatorname{lin}}(\{P_n:n\neq m\}\cup\{Q_1,Q_2\})$ so 
\[A\neq\overline{\operatorname{lin}}(\{P_n:n\neq m\}\cup\{Q_1,Q_2\}).\]
Thus, $A$ is not maximal, a contradiction.

If $P\in L(s',s)$ is a nonzero projection orthogonal to all $P_n$, then, using similar arguments as above, we get that $\overline{\operatorname{lin}}(\{P_n\}_{n\in\N}\cup\{P\})$ is a closed 
commutative ${}^*$-subalgebra of $L(s',s)$ bigger than $A$ which leads to a contradiction.



(ii)$\Rightarrow$(iii). Since $A$ is commutative thus $A\subset A'$. Now, let us suppose that there is $x\in A'\setminus A$. By Proposition \ref{prop_Z'},
$x^*\in A'$ so $x+x^*,x-x^*\in A'$, and, moreover, $x^*\notin A$. Since $x=\frac{1}{2}(x+x^*)+\frac{1}{2}(x-x^*)$ thus $x+x^*\notin A$ or 
$x-x^*\notin A$. Without loss of 
generality let us assume that $x+x^*\notin A$. Operator $x+x^*$ is self-adjoint hence it has the spectral representation of the form 
$\sum_{m=1}^\infty\mu_mQ_m$.
Then, by Propositions \ref{prop_x_in_A_iff_P_n_in_A} and \ref{prop_Z'}, 
$Q_m\in A'$ for all $m\in\N$ and there exists $m_0$ for which $Q_{m_0}\notin A$ (otherwise $x+x^*\in A$).
Let $J:=\{n:P_n\preceq Q_{m_0}\}$ (see the definition of $\preceq$ in Lemma \ref{lem_minimal_projections}). Since $Q_{m_0}$ is finite dimensional thus
$J$ is finite. One can easily check that $Q_{m_0}-\sum_{j\in J}P_j$ is a projection (if $J=\emptyset$, then $\sum_{j\in J}P_j:=0$). Moreover, 
\begin{equation}\label{eq_Q_m_0}
\Big(Q_{m_0}-\sum_{j\in J}P_j\Big)P_k=0
\end{equation}
for all $k\in\N$.
Indeed, if $k\in J$, then from the definition of the relation $\preceq$, $Q_{m_0}P_k=P_k$, so 
\[\Big(Q_{m_0}-\sum_{j\in J}P_j\Big)P_k=Q_{m_0}P_k-P_k=0.\]
Let $k\notin J$. We have $Q_{m_0}P_k=P_kQ_{m_0}$ because $Q_{m_{0}}\in A'$. This implies that 
$Q_{m_0}P_k$ is a projection and $\im Q_{m_0}P_k=\im Q_{m_0}\cap \im P_k$. 
Therefore, since $P_k$ are one-dimensional thus $Q_{m_0}P_k=P_k$ or $Q_{m_0}P_k=0$. But, from the assumption, $Q_{m_0}P_k\neq P_k$ so $Q_{m_0}P_k=0$. Now,
\[\Big(Q_{m_0}-\sum_{j\in J}P_j\Big)P_k=Q_{m_0}P_k=0.\]

Since the sequence $(P_n)_{n\in\N}$ is complete, (\ref{eq_Q_m_0}) implies that $Q_{m_0}-\sum_{j\in J}P_j=0$. Hence $Q_{m_0}\in A$, a contradiction.

(iii)$\Rightarrow$(i). Follows directly from the definition of the comutant of $A$.
\qed

\begin{Rem}\label{rem1}
(i) Since $(P_n)_{n\in\N}$ is a sequence of pairwise orthogonal one-dimensional projections so $P_n=\langle\cdot,e_n\rangle e_n$, 
where $(e_n)_{n\in\N}\subset s$ is an orthonormal system in $l_2$. Then $\lambda^\infty(||P_n||_q)=\lambda^\infty(|e_n|_q)$ as Fr\'echet ${}^*$-algebras.
Indeed, from the H\"older inequality, if $\xi\in s$, $q\in\N_0$, then 
\[|\xi|_q^2\leq|\xi|_{l_2}|\xi|_{2q}.\]
Hence
\[1\leq|e_n|_q\leq|e_n|_q^2=||P_n||_q=|e_n|_q^2\leq|e_n|_{2q}.\]
This implies that $\lambda^\infty(||P_n||_q)=\lambda^\infty(|e_n|_q)$ as Fr\'echet spaces, and, since the algebra operations are the same in both algebras, 
$\lambda^\infty(||P_n||_q)=\lambda^\infty(|e_n|_q)$ 
as Fr\'echet ${}^*$-algebras.

(ii) The sequence $(e_n)_{n\in\N}$ from the previous item need not to be an orthonormal basis of $l_2$. Indeed, let $(e_n)_{n\in\N}$ be an orthonormal
basis of $l_2$ such that $e_n\in s$ for $n\in\N\setminus\{1\}$ and $e_1\notin s$. Clearly, $(e_n)_{n\in\N\setminus\{1\}}$ 
is not an orthonormal basis of $l_2$ and $(\langle\cdot,e_n\rangle e_n)_{n\in\N\setminus\{1\}}$ is a complete sequence in $L(s',s)$.

(iii) By the Kuratowski-Zorn lemma, every closed commutative ${}^*$-subalgebra of $L(s',s)$ is contained in some maximal commutative 
${}^*$-subalgebra of $L(s's)$. If $\{P_n\}_{n\in\N}$ is a sequence of pairwise orthogonal finite dimensional projections, then, by Proposition 
\ref{prop_alg=lin}, $\alg(\{P_n\}_{n\in\N})$ is a closed commutative ${}^*$-subalgebra of $L(s',s)$ so it is contained in some maximal commutative
${}^*$-subalgebra $\alg(\{Q_n\}_{n\in\N})$ of $L(s's)$, where $\{Q_n\}_{n\in\N}$ is a complete sequence of one-dimensional 
projections in $L(s',s)$ (see Theorems \ref{th_commutative_subalg} and \ref{th_max_com}). Now, applying Lemma \ref{lem_minimal_projections} (i), it is
easy to show that the sequence $\{P_n\}_{n\in\N}$ can be extended to some complete sequence of projections belonging to $L(s',s)$. 

\end{Rem}

\begin{Cor}
Let $A$ be one of the following Fr\'echet ${}^*$-algebras with the pointwise multiplication:
\begin{enumerate}
\item[(i)] the algebra $\mathcal S(\R^n)$ of rapidly decreasing smooth functions;
\item[(ii)] the algebra $\mathcal D(K)$ of test functions with support in a compact set $K\subset\R^n$, 
$\operatorname{int}(K)\neq\emptyset$;
\item[(iii)] the algebra $C^\infty_a(M)$ of smooth functions on a compact smooth manifold $M$ vanishing at $a\in M$;
\item[(iv)] the algebra $C^\infty_a(\overline{\Omega})$ of smooth functions on $\overline{\Omega}$ vanishing at $a\in\Omega$, where $\Omega\neq\emptyset$ is an open bounded subset of $\R^n$ with $C^1$-boundary;
\item[(v)] the algebra $\mathcal E_a(K)$ of Withney jets on a compact set $K\subset\R^n$ with the extension property, flat at $a\in K$ and such that 
$\operatorname{int}(K)\neq\emptyset$. 
\end{enumerate}
Then $A$ is isomorphic to $s$ as Fr\'echet space but it is not isomorphic to a closed commutative ${}^*$-subalgebra of $L(s',s)$ as Fr\'echet 
${}^*$-algebra.
\end{Cor}
\proof
It is well known that spaces from items (i)-(v) are isomorphic to the space $s$ as Fr\'echet space (see e.g. \cite[Ch. 31]{MeV}, \cite[Satz 4.1]{V2}). 

To prove the second assertion let us compare sets of ${}^*$-multiplicative functionals. 
If $A$ is one of the spaces from items (i)-(v), then every point evaluation functional on $A$ 
is ${}^*$-multiplicative and since the cardinal number of the underlying space is continuum, thus the cardinal number of the set of 
${}^*$-multiplicative functionals 
on $A$ is not less than continuum. On the other hand, by Corollary \ref{cor_coefficient_functionals}, the set of ${}^*$-multiplicative functional
on any infinite dimensional closed commutative ${}^*$-subalgebra of $L(s',s)$ is at most countable, hence no of the spaces from items (i)-(v) is isomorphic to $A$.
\qed

It is clear that the algebra $s$ with the pointwise multiplication is a ${}^*$-subalgebra of $L(s',s)$ (consider, for example, diagonals operators). 
The previous corollary shows that it is not the case for many other interesting Fr\'echet ${}^*$-algebras isomorphic to $s$ (as Fr\'echet spaces).
This leads to the following: 
\begin{OPr}\label{opr1}
Does every closed commutative ${}^*$-subalgebra of $L(s',s)$ is isomorphic to some closed ${}^*$-subalgebra of the 
algebra $s$ with the pointwise multiplication? 
\end{OPr}

We say that an orthonormal system $(e_n)_{n\in\N}$ in $l_2$ is \emph{$s$-complete}, if every $e_n$ belongs to $s$ and 
for every $\xi\in s$ the following implication holds:
if $\langle \xi,e_n\rangle=0$ for every $n\in\N$, then $\xi=0$. One can easily show that an orthonormal system $(e_n)_{n\in\N}$ is $s$-complete 
if and only if the sequence of projections $(\langle\cdot,e_n\rangle e_n)_{n\in\N}$ is complete in $L(s',s)$.
Hence, by Theorem \ref{th_max_com} and Remark \ref{rem1}, Problem \ref{opr1} is equivalent to the following.
\begin{OPr}
Let $(e_n)_{n\in\N}$ be $s$-complete, orthonormal system in $l_2$. 
Is the algebra $\lambda^\infty(|e_n|_q)$ isomorphic to some closed ${}^*$-subalgebra of the algebra $s$?
 \end{OPr}


\section{Functional calculus in $L(s',s)$}


If $x$ is a normal operator in $L(s',s)\subset K(l_2)$ and $f$ is a continuous function on the spectrum $\sigma(x)$ of $x$ vanishing at zero, then,
from the continuous functional calculus for normal operators (see e.g. \cite[Prop. 17.20]{MeV}), $f(x)$ is the uniquely determined operator in $K(l_2)$.
In this section, we want to describe these functions $f$ for which $f(x)$ is again in $L(s',s)$.  

From the general theory of Fr\'echet $m$-locally convex algebras we get the holomorphic functional calculus on $L(s',s)$ 
(see Prop. \ref{prop_6} and \cite[Lemma 1.3]{Phill}, \cite[Th. 12.16]{Zel}). 
Precisely, if $x$ is an arbitrary operator in $L(s,s)$ and $f$ is a holomorphic function on an open neighborhood $U$ of $\sigma(x)$ with $f(0)=0$, then
$f(x)\in L(s',s)$, and, moreover, the map $\Phi\colon H_0(U)\to L(s',s)$, $f\mapsto f(x)$ is a continuous homomorphism ($H_0(U)$ stands for 
the space of holomorphic functions vanishing at zero). 

Let us recall that $f\colon X\to \C$ ($X\subset \C$, $0\in X$) is a \emph{H\"older continuous at zero function} if there is $\theta\in(0,1]$ and $C>0$ 
such that $|f(t)-f(0)|\leq C|t|^\theta$ for all $t$ in the neighborhood of 0. As a consequence of Theorem \ref{th_spectral} 
we get the following H\"older continuous functional calculus for normal operators in $L(s',s)$.

\begin{Cor}\label{cor_holder}
If $x\in L(s',s)\subset K(l_2)$ is normal, then for every H\"older continuous at zero function $f\colon\sigma(x)\to\C$ with $f(0)=0$, 
the operator $f(x)\in L(s',s)$ as well. 
In particular, for every normal operator $x\in L(s',s)$ with $\sigma(x)\subset [0,\infty)$ and $\theta\in (0,\infty)$, $x^\theta\in L(s',s)$.  
\end{Cor}
\proof
Let $x=\sum_{n\in\cN}\lambda_nP_n$ be a normal operator in $L(s',s)$ with the nonnegative spectrum and let $\theta\in(0,\infty)$. If $\theta\in(0,1]$, then,
by Theorem \ref{th_spectral}, $x^\theta=\sum_{n=1}^\infty\lambda_n^\theta P_n\in L(s',s)$ and for $\theta\in(1,\infty)$ we have
$x^\theta=x^{\lfloor\theta\rfloor}\cdot x^{\theta-\lfloor\theta\rfloor}\in L(s',s)$, where $\lfloor\theta\rfloor$ is the floor of $\theta$.

Now, let $x=\sum_{n\in\cN}\lambda_nP_n\in L(s',s)$ be normal and let $f\colon\sigma(x)\to\C$ be a H\"older continuous at zero function with $f(0)=0$. 
Then $|f|\leq C|\cdot|^\theta$ for some $C>0$ and $\theta\in(0,1]$. Hence
$\sum_{n\in\cN} ||f(\lambda_n)P_n||_q\leq C\sum_{n\in\cN}|\lambda_n|^\theta ||P_n||_q<\infty$ so, by Corollary \ref{cor_characterization_of_L(s',s)}, $f(x)\in L(s',s)$.
\qed

For a normal operator $x$ in $L(s',s)$ with the spectral representation $x=\sum_{n=1}^\infty\lambda_nP_n$, we define the function space
\[C_s(\sigma(x)):=\{f:\sigma(x)\To\C: f(0)=0, (f(\lambda_n))_{n\in\N}\in\lambda^\infty(||P_n||_q)\}.\]
It is easy to show that the space $C_s(\sigma(x))$ with a system $(c_q)_{q\in\N_0}$, $c_q(f):=\sup{|f(\lambda_n)|||P_n||_q}$ of seminorms,
pointwise multiplication and conjugation is a Fr\'echet ${}^*$-algebra.

\begin{Th}\label{th_functional_calculus}
If $x=\sum_{n=1}^\infty\lambda_nP_n$ is an infinite dimensional normal operator in $L(s',s)$, then the map 
\[\Phi:C_s(\sigma(x))\To\alg(x),\quad \Phi(f):=f(x):=\sum_{n=1}^\infty f(\lambda_n)P_n\]
is a Fr\'echet algebra isomorphism such that
$\Phi(\id)=x$ and $\Phi(\overline{f})=\Phi(f)^*$.
\end{Th}
\proof
By Theorem \ref{th_commutative_subalg}, $\Phi$ is well defined, and, of course, $\Phi(\id)=x$, $\Phi(\overline{f})=\Phi(f)^*$.
The space $\alg(x)$ is a nuclear Fr\'echet space (as a closed subspace of nuclear Fr\'echet space $L(s',s)$) 
so $\lambda^\infty(||P_n||_q)\cong\alg(x)$ (see Theorem \ref{th_commutative_subalg}) is a nuclear Fr\'echet space as well.
Thus, by the Grothendieck-Pietsch theorem (see e.g. \cite[Th. 28.15]{MeV}), for given $q\in\N_0$ one can 
find $r\in\N_0$ such that $C:=\sum_{n=1}^\infty\frac{||P_n||_q}{||P_n||_r}<\infty$. Hence
\begin{align*}
||\Phi(f)||_q&\leq\sum_{n=1}^\infty|f(\lambda_n)|||P_n||_q=\sum_{n=1}^\infty|f(\lambda_n)|||P_n||_r\frac{||P_n||_q}{||P_n||_r}\\
&\leq \sup_{n\in\N}|f(\lambda_n)|||P_n||_r\cdot\sum_{n=1}^\infty\frac{||P_n||_q}{||P_n||_r}=Cc_r(f),
\end{align*}
thus $\Phi$ is continuous. 

Clearly, $\Phi$ is injective. To prove that $\Phi$ is also surjective, let us take $y\in\alg(x)$. 
Since, by Theorem \ref{th_commutative_subalg}, $(P_n)_{n\in\N}$ is
a Schauder basis, so there is a sequence
$(\mu_n)_{n\in\N}$ such that $y=\sum_{n=1}^\infty\mu_nP_n$. Let $g(\lambda_n):=\mu_n$ for $n\in\N$. Then 
\[\sup_{n\in\N}|g(\lambda_n)|||P_n||_q=\sup_{n\in\N}|\mu_n|||P_n||_q<\infty,\]
hence $g\in C_s(\sigma(x))$ and of course $\Phi(g)=y$.
\qed

\flushleft{\textbf{Acknowledgements}. The author is very indebted to Pawe\l{} Doma\'nski and Krzysztof Piszczek for several helpful conversations 
and especially for discussions on their papers \cite{Dom}, \cite{Piszczek}.

\vspace{1cm}
\begin{minipage}{7.5cm}
T. Cia\'s

Faculty of Mathematics and Comp. Sci.

A. Mickiewicz University in Pozna{\'n}

Umultowska 87

61-614 Pozna{\'n}, POLAND

e-mail: tcias@amu.edu.pl
\end{minipage}\


{\small
\begin{thebibliography}{99}

\bibitem{BhIO} S. J. Bhatt, A. Inoue, H. Ogi, 
\emph{Spectral invariance, $K$-theory isomorphism and an application to the differential structure of $C^*$-algebras}.
J. Operator Theory \textbf{49} (2003), no. 2, 389-405.

\bibitem{BlCun} B. Blackadar, J. Cuntz,
\emph{Differential Banach algebra norms and smooth subalgebras of $C^*$-algebras}.
J. Operator Theory \textbf{26} (1991), no. 2, 255-282.

\bibitem{Con} A. Connes,
\emph{Noncommutative geometry}.
Academic Press, Inc., San Diego, CA, 1994.

\bibitem{Conw} J. B. Conway,
\emph{A course in functional analysis}.
Second edition. Graduate Texts in Mathematics, \textbf{96}. Springer-Verlag, New York, 1990.

\bibitem{Cuntz} J. Cuntz,
\emph{Bivariante K-Theorie f{\"u}r lokalkonvexe Algebren und der Chern-Connes-Charakter}.
Doc. Math. \textbf{2} (1997), 139-182.
 
\bibitem{Dom} P. Doma\'nski,
\emph{Algebra of smooth operators}.
Unpublished note.

\bibitem{ElNatNest} G. A. Elliot, T. Natsume, R. Nest,
\emph{Cyclic cohomology for one-parameter smooth crossed products}.
Act. Math. \textbf{160} (1998), 285-305.



\bibitem{Fra} M. Fragoulopoulou,
\emph{Topological algebras with involution}.
North-Holland Mathematics Studies, \textbf{200}. Elsevier Science B.V., Amsterdam, 2005.


\bibitem{GlockLkamp} H. Gl{\"o}ckner, B. Langkamp,
\emph{Topological algebras of rapidly decreasing matrices and generalizations}.
Topology Appl. \textbf{159} (2012), no. 9, 2420-2422. 

\bibitem{Kot} G. K\"othe,
\emph{Topological vector spaces II}.
Springer-Verlag, Berlin-Heidelberg-New York 1979.

\bibitem{MeV} R. Meise, D. Vogt,
\emph{Introduction to functional analysis}.
Oxford University Press, New York 1997.

\bibitem{Phill} N. C. Phillips,
\emph{K-theory for Fr\'echet algebras}.
Internat. J. Math. \textbf{2} (1991), no. 1, 77-129. 

\bibitem{Piszczek} K. Piszczek,
\emph{On a property of $PLS$-spaces inherited by their tensor products}.
Bull. Belg. Math. Soc. Simon Stevin \textbf{17} (2010), 155-170.

\bibitem{Schw} L. B. Schweitzer,
\emph{A short proof that $M_n(A)$ is local if $A$ is local and Fr\'echet}.
Internat. J. Math. \textbf{3} (1992), no. 4, 581-589.
 
\bibitem{V2} D. Vogt,
\emph{Ein Isomorphiesatz f\"ur Potenzreihenr\"aume}.
Arch. Math. (Basel) \textbf{38} (1982), no. 6, 540-548. 

\bibitem{V1} D. Vogt,
\emph{On the functors $\operatorname{Ext}^1(E,F)$ for Fr\'echet spaces}.
Studia Math. \textbf{85} (1987), no. 2, 163-197. 

\bibitem{Zel} W. \.Zelazko,
\emph{Selected topics in topological algebras}.
Lectures 1969/1970, Lectures Notes Series, No. \textbf{31}. Matematisk Institut, Aarhus Universitet, Aarhus 1971. 

\end{thebibliography}
}
\end{document}